\documentclass[12pt,twoside]{article}

\usepackage{amssymb}
\usepackage{isolatin1}
\usepackage{epsfig}

\usepackage{graphicx}

\usepackage{color}

\newtheorem{thm}{Theorem}
\newtheorem{definition}[thm]{Definition}       
\newtheorem{lem}[thm]{Lemma}
\newtheorem{corollary}[thm]{Corollary}
\newtheorem{example}[thm]{Example}

\newtheorem{proposition}[thm]{Proposition}

\newtheorem{fremdersatz}{Theorem}

\newenvironment{satOrig}[1]
        {\pagebreak[2] \begin{fremdersatz} {\bf #1} \quad\sl}
        { \end{fremdersatz}}

\newenvironment{proofof}[1]
        {\pagebreak[2] \vspace{-1pt}{\bf Proof#1.}  }
        {\hfill $\blacksquare$ \vspace{2pt}}

\def\parmod{\parskip=2pt plus1pt minus1pt}

\newenvironment{proof}
        {\pagebreak[2] \vspace{-1pt}{\bf Proof.}  }
        {\hfill $\blacksquare$ \vspace{2pt}}

\def\fa{\mathcal{F}}

\def\nat{{\rm I\! N}}

\def\co{{\mathbb C}}
\def\re{{\mathbb R}}

\def\l{\left}
\def\r{\right}
\def\gl{\left\{}
\def\gr{\right\}}
\def\kl{\left(}
\def\kr{\right)}

\def\limj{\lim_{j\to\infty}}

\def\limn{\lim_{n\to\infty}}

\def\abb{\longrightarrow}

\renewcommand{\rho}{\varrho}
\renewcommand{\phi}{\varphi}
\renewcommand{\epsilon}{\varepsilon}

\def\beq{\begin{equation}}
\def\eeq{\end{equation}}
\def\beqar{\begin{eqnarray}}
\def\eeqar{\end{eqnarray}}
\def\beqaro{\begin{eqnarray*}}
\def\eeqaro{\end{eqnarray*}}
\def\bsat{\begin{thm}}
\def\esat{\end{thm}}
\newcommand{\bsatorig}{\begin{satOrig}}
\newcommand{\esatorig}{\end{satOrig}}
\def\blem{\begin{lem}}
\def\elem{\end{lem}}
\def\bkor{\begin{corollary}}
\def\ekor{\end{corollary}}
\def\bprop{\begin{proposition}}
\def\eprop{\end{proposition}}
\def\bdefin{\begin{definition}}
\def\edefin{\end{definition}}
\def\bbew{\begin{proof}}
\def\ebew{\end{proof}}
\def\bbewo{\begin{proofof}}
\def\ebewo{\end{proofof}}
\def\bgar{\begin{array}}
\def\ear{\end{array}}

\def\bex{\begin{example}}
\def\eex{\end{example}}

\renewcommand{\rho}{\varrho}
\renewcommand{\phi}{\varphi}

\begin{document}

\thispagestyle{plain}

\begin{center}
{\LARGE\bf On the growth of real functions and their derivatives \\[25pt] } 

{\Large \it Jürgen Grahl and Shahar Nevo}
\end{center}

\begin{abstract}

\vspace{-12pt}

We show that for any $k$-times continuously differentiable function
 $f:[a,\infty)\abb\re$, any integer $q\ge 0$  and any $\alpha>1$ the
 inequality
$$\liminf_{x\to\infty} \frac{x^k \cdot\log x\cdot \log_2 x\cdot\dots\cdot
  \log_q x \cdot f^{(k)}(x)}{1+|f(x)|^\alpha}\le 0 $$
holds.
\end{abstract}

{\it 2000 Mathematics Subject Classification: 26D10}

In \cite{GrahlNevo-Spherical}, \cite{LiuNevoPang},  
 \cite{ChenNevoPang}, \cite{GNP-NonExplicit},
\cite{GN-Marty} and \cite{BarGrahlNevo} we had studied differential
inequalities in the context of complex analysis, more precisely with
respect to the question whether they constitute normality (or at
least quasi-normality) in the sense of Montel. 

\bsatorig{\cite{ChenNevoPang}} \label{CNP}
Let $\alpha> 1$ and $C>0$ be real numbers and $k\ge 1$ be an
integer. Let $\fa$ be a family of meromorphic functions in some domain
$D$ in the complex plane such that  
\beq\label{LowerEstimate}
\frac{|f^{(k)}|}{1+|f|^\alpha}(z)\ge C \qquad \mbox{ for all }
z\in D \mbox{ and all } f\in\fa.
\eeq
Then $\fa$ is normal.
\esatorig

This result doesn't hold any longer if $\alpha>1$ is replaced by
$\alpha=1$ as easy examples demonstrate. However, at least for $k=1$ 
condition (\ref{LowerEstimate}) implies quasi-normality if
$\alpha=1$ \cite{LiuNevoPang}. Furthermore, in \cite{BarGrahlNevo} we
had shown that the condition  
\beq\label{LowerBound}
\frac{|f^{(k)}|}{1+|f^{(j)}|^\alpha}(z)\ge C \qquad \mbox{ for all }
z\in D
\eeq
(where $k>j\ge 0$ are integers, $\alpha>1$ and $C>0$) implies
quasi-normality.  

As to entire functions (i.e. functions analytic in the whole complex
plane), it is almost obvious that they cannot satisfy a differential
inequality like (\ref{LowerEstimate}). Indeed, if $f$ is entire and 
$|f^{(k)}|(z)\ge C\cdot (1+|f(z)|^\alpha)$ for all $z\in \co$, then in
particular $|f^{(k)}(z)|\ge C$ for all $z\in\co$,  so $f^{(k)}$ is
constant by Picard's (or Liouville's) theorem. But then $f$ is a non-constant
polynomial, and one obtains a contradiction for $z\to\infty$ provided
that $\alpha >0$. 
 
These considerations motivated us to look at the differential
inequality (\ref{LowerEstimate}) in the context of {\bf real} analysis, a
problem that doesn't seem to have been studied so far.  For
real-valued functions on unbounded intervals we have the following
result which turns out to be sharp in a certain sense. Here, $\log_p
x$ denotes the $p$-times iterated natural logarithm, defined
recursively by $\log_0 x:=x$ and $\log_p x:=\log(\log_{p-1} x)$ for
$p\ge 1$.

\bsat{} \label{mainresult}
Let $k\ge 1$ and $q\ge 0$ be integers, $\alpha>1$, $a\in\re$ and
$f:[a,\infty)\abb\re$ a $k$-times continuously differentiable function. Then
\beq\label{MainAssertion1}
\liminf_{x\to\infty} \frac{x^k \cdot\log x\cdot \log_2 x\cdot\dots\cdot
  \log_q x \cdot f^{(k)}(x)}{1+|f(x)|^\alpha}\le 0 
\eeq
and
\beq\label{MainAssertion2}
\liminf_{x\to\infty} \frac{x^k\cdot\log x\cdot \log_2 x\cdot\dots\cdot
  \log_q x \cdot |f^{(k)}(x)|}{1+|f(x)|^\alpha}= 0. 
\eeq
\esat

\bex
This result is best possible in the sense that it is not longer valid
if $\log_q x$ is replaced by $(\log_q x)^\beta$ with any
$\beta>1$. This can be seen by considering the function
$f:[a,\infty[\abb\re$ defined by
$$f(x):=(-1)^{k-1}\cdot\int_{a}^{x}\int_{x_k}^{\infty}\dots\int_{x_2}^{\infty}
\frac{1}{x_1^k\cdot \log x_1\cdot\ldots\cdot \log_q
  x_1}\;dx_1\,\dots\,dx_k,$$
where $a>0$ is chosen sufficiently large\footnote{Another, related example is
  $f(x):=\log_{q+1} x$. However, it is more difficult to verify that
  it has the desired properties than for the example given above.}. Indeed, for $x\ge a$ we have
\beqaro
|f(x)|&\le& 
\int_{a}^{x}\frac{1}{\log x_k\cdot\ldots\cdot \log_q x_k}\kl\int_{x_k}^{\infty}\dots\int_{x_2}^{\infty}
\frac{1}{x_1^k}\;dx_1\,\dots\,dx_{k-1}\kr\;dx_k\\
&=& \frac{1}{(k-1)!}\int_{a}^{x}\frac{1}{\log x_k\cdot\ldots\cdot
  \log_q x_k}\cdot \frac{1}{x_k}\;dx_k\\
&=& \frac{1}{(k-1)!}\cdot \log_{q+1} x
\eeqaro
and of course
$$f^{(k)}(x)=\frac{1}{x^k\cdot \log x\cdot\ldots\cdot \log_q x},$$
hence for any $\alpha,\beta>1$ 
$$\frac{x^k \cdot\log x\cdot \log_2 x\cdot\dots\cdot (\log_q x)^\beta \cdot
  f^{(k)}(x)}{1+|f(x)|^\alpha}
\ge \frac{(\log_q x)^{\beta-1}}{1+\kl\frac{1}{(k-1)!}\cdot \log_{q+1} x\kr^\alpha}
\to \infty \qquad (x\to\infty).$$
So (\ref{MainAssertion1}) does not hold, and neither does (\ref{MainAssertion2}).
\eex

Of course, the appearance of the terms $\log x\cdot \log_2
x\cdot\dots\cdot \log_q x$ in Theorem \ref{mainresult}, where $\log_q
x$ cannot be replaced by $(\log_q x)^\beta$ with $\beta>1$, is
reminescent of the well-known fact from basic calculus that for any
natural number $q$ the infinite series $\sum_{k=k_0}^{\infty} (k\log
k\cdot\dots\cdot \log_{q-1} k\cdot (\log_q k)^\beta)^{-1}$ (where
$k_0$ is chosen sufficiently large) is convergent for $\beta>1$ and
divergent for $0<\beta\le 1$ and that a corresponding result holds for
the improper integral $\int_{x_0}^{\infty} (x\cdot \log
x\cdot\dots\cdot \log_{q-1} x\cdot (\log_q x)^\beta)^{-1}\;dx$. This
resemblance seems to be more than coincidence as Case 3 of the proof
of (\ref{MainAssertion1}) reveals: It makes crucial use of the
divergence of $\int_{x_0}^{\infty} (x\cdot \log x\cdot\dots\cdot
\log_q x)^{-1}\;dx$.

\bbew{}
Our main efforts are required to prove (\ref{MainAssertion1}). Then
(\ref{MainAssertion2}) will be an easy consequence from
(\ref{MainAssertion1}).   

We want to prove (\ref{MainAssertion1}) by induction w.r.t. $q$. However, the start of
our induction is to consider $\frac{f^{(k)}(x)}{1+|f(x)|^\alpha}$
rather than $\frac{x^k\cdot f^{(k)}(x)}{1+|f(x)|^\alpha}$ (which would
be the case $q=0$). So we have to introduce a unifying notation
first. For given $k\ge 1$, we set 
$$P_{-1}(x):=1 \qquad \mbox{ and } \qquad 
P_q(x):=x^k\cdot \prod_{j=1}^{q} \log_j x \quad\mbox{ for } q\ge 0.$$
In particular, $P_0(x)=x^k$. Then (\ref{MainAssertion1}) has the form 
$$\liminf_{x\to\infty} \frac{P_q(x)\cdot
  f^{(k)}(x)}{1+|f(x)|^\alpha}\le0.$$

First we consider the case $q=-1$. Let's assume the assertion is
wrong. Then there is an $\varepsilon>0$ 
and an $a_0\ge 0$ such that 
$$f^{(k)}(x)\ge \varepsilon\cdot \kl 1+|f(x)|^\alpha\kr \qquad \mbox{
  for all } x\ge a_0.$$
From $f^{(k)}(x)\ge \varepsilon$ for all $x\ge a_0$ one easily sees
  that there is some $x_1\ge a_0$ such that $f^{(k)}(x)>0, 
  f^{(k-1)}(x)>0,\dots, f'(x)>0, f(x)>0$ for all $x\ge x_1$. In
  particular, $f$ is strictly increasing (i.e. one-to-one) on
  $[x_1,\infty[$ and  $\lim_{x\to\infty} f(x)=\infty$. We choose a
  natural number $n$ such that $(\alpha-1)\cdot n>k-1$. Then there is
  a natural number $j_0$ such that $f([x_1,\infty[)$ contains the
  interval $[j_0^n,\infty[$. For $j\ge j_0$ we set 
$$r_j:=f^{-1}(j^n).$$
Then $(r_j)_j$ is strictly increasing and unbounded, and by the mean
value theorem, applied to $\varphi(t):=t^n$, we have  
$$f(r_{j+1})-f(r_j)=(j+1)^n-j^n\le n \cdot (j+1)^{n-1} 
\qquad \mbox{ for all } j\ge j_0.$$ 
On the other hand, for $j\ge j_0$ we deduce from the fundamental
theorem of calculus
\beqaro
f(r_{j+1})-f(r_j)&=& \int_{r_j}^{r_{j+1}} f'(x_1)\;dx_1 \\
&= & \int_{r_j}^{r_{j+1}} \kl f'(r_j)+\int_{r_j}^{x_1} f''(x_2)\;dx_2 \kr\;dx_1\\
&\ge & \int_{r_j}^{r_{j+1}} \int_{r_j}^{x_1} f''(x_2)\;dx_2
\;dx_1\\[9pt]
&\ge & \dots \\[9pt]
&\ge & \int_{r_j}^{r_{j+1}} \int_{r_j}^{x_1}\dots\int_{r_j}^{x_{k-1}}
f^{(k)}(x_k)\;dx_k\dots dx_2 dx_1\\
&\ge & \varepsilon\cdot \int_{r_j}^{r_{j+1}} \int_{r_j}^{x_1}\dots\int_{r_j}^{x_{k-1}}
\kl 1+ f^\alpha(x_k)\kr\;dx_k\dots dx_2 dx_1\\
&\ge & \varepsilon\cdot \int_{r_j}^{r_{j+1}} \int_{r_j}^{x_1}\dots\int_{r_j}^{x_{k-1}}
f^\alpha(r_j)\;dx_k\dots dx_2 dx_1\\
&=& \varepsilon\cdot j^{\alpha n}\cdot\frac{1}{k!}\cdot (r_{j+1}-r_j)^k.
\eeqaro
Combining these two estimates yields 
$$n \cdot (j+1)^{n-1} \ge \frac{\varepsilon}{k!}\cdot j^{\alpha
  n}\cdot (r_{j+1}-r_j)^k,$$
hence
$$r_{j+1}-r_j 
\le \kl \frac{n\cdot k!}{\varepsilon}\cdot\frac{(j+1)^{n-1}}{j^{\alpha  n}}\kr^{1/k}
\le \kl \frac{n\cdot k!\cdot 2^{n-1}}{\varepsilon}\kr^{1/k}\cdot
j^{-((\alpha-1)\cdot n+1)/k}.$$
Here, by our choice of $n$, $((\alpha-1)\cdot n+1)/k>1$, so the series
$\sum_{j=j_0}^{\infty} j^{-((\alpha-1)\cdot n+1)/k}$ converges. Hence also the
telescope series $\sum_{j=j_0}^{\infty} (r_{j+1}-r_j)=\limj
r_j-r_{j_0} $ converges, contradicting $\limj r_j=\infty$. This proves
(\ref{MainAssertion1}) for $q=-1$. 

Now let some $q\ge 0$ be given and assume that (\ref{MainAssertion1}) is
true for $q-1$ instead of $q$ and for all $k$-times differentiable
functions $f:[0,\infty)\abb\re$. We assume there is a $k$-times differentiable
function $f:[0,\infty)\abb\re$ and an $\varepsilon>0$ such that
\beq\label{ContradAssumpt}
P_q(x)\cdot{ f^{(k)}}(x)\ge \varepsilon\cdot \kl1+|f(x)|^\alpha\kr
\eeq
holds for all $x$ large enough. Then in particular $f^{(k)}(x)>0$ for
all large enough $x$, so $f^{(k-1)}$ is increasing, and we easily see
by induction that $f^{(k-1)}, f^{(k-2)},\dots,f',f$ are strictly monotonic on
an appropriate interval $[x_0,\infty)$ (for large enough $x_0$). So the
limits
$$L_j:=\lim_{x\to\infty} f^{(j)}(x) \qquad (j=0,\dots,k-1)$$
exist. (They might be $+\infty$ or $-\infty$.)

In the following we will apply the induction hypothesis to the 
function
$$g(t):=f(e^t)$$
and will use that 
\beq\label{DerivG}
g^{(k)}(t)=f^{(k)}(e^t)\cdot e^{kt}+\sum_{j=1}^{k-1} c_j
f^{(j)}(e^t)\cdot e^{jt}
\eeq
for certain constants $c_j\ge 0$. (This is easily seen by induction.)

By the mean value theorem, for all $n\in\nat$ there is a
$\zeta_n\in[n,2n]$ such that 
\beq\label{MVT}
n\cdot |f^{(k)}(\zeta_n)|=|f^{(k-1)}(2n)-f^{(k-1)}(n)|.
\eeq
Here of course we have $\limn \zeta_n=\infty$. 

Now we consider several cases. 

{\bf Case 1: } $L_{k-1}\ne 0$.

Since $f^{(k-1)}$ is increasing, we either have $L_{k-1}\in\re$ or
$L_{k-1}=+\infty$. 

{\bf Case 1.1: } $L_{k-1}\in\re$, w.l.o.g. $L_{k-1}>0$.

Then we have
$$\frac{1}{2}\cdot L_{k-1}\le f^{(k-1)}(x) \le 2 L_{k-1} \qquad \mbox{ for large enough } x,$$
hence
$$\frac{1}{3(k-1)!}\cdot L_{k-1}\cdot x^{k-1}\le f(x) \le
\frac{3}{(k-1)!} L_{k-1}\cdot x^{k-1} \qquad \mbox{ for large enough }
x.$$
Using the lower estimate, we conclude that for large enough $x$ 
\beq\label{Case1.1a}
0\le P_q(x)\cdot \frac{1}{x}\cdot \frac{1}{1+|f(x)|^\alpha}
\le \frac{x^{(k-1)(1+\alpha)/2}}{1+|f(x)|^\alpha} \abb 0 \quad
(x\to\infty).
\eeq
(Here it is crucial that $1<\frac{1}{2}\cdot (1+\alpha)<\alpha$.)
Furthermore, 
\beq\label{Case1.1b}
0\le \zeta_n\cdot |f^{(k)}(\zeta_n)|
\le 2n\cdot |f^{(k)}(\zeta_n)|
=2\cdot |f^{(k-1)}(2n)-f^{(k-1)}(n)|\abb 0 \quad (n\to\infty)
\eeq
since $L_{k-1}$ is finite. Multiplying (\ref{Case1.1a}) and
(\ref{Case1.1b}) gives 
$$0\le P_q(\zeta_n)\cdot
\frac{|f^{(k)}(\zeta_n)|}{1+|f(\zeta_n)|^\alpha}\abb0 \quad
(n\to\infty).$$
This is a contradiction to (\ref{ContradAssumpt}). 

{\bf Case 1.2: } $L_{k-1}=+\infty$.

Then for large enough $x$ we have $f^{(k-1)}(x)\ge 1, f^{(k-2)}(x)\ge
1, \dots, f'(x)\ge 1, f(x)\ge 1$ (and
$L_{k-2}=\dots=L_1=L_0=+\infty$). By applying the induction hypothesis
to $g$, using (\ref{DerivG}) and substituting $t=\log x$ we obtain
\beqaro
0&\ge & \liminf_{t\to+\infty} P_{q-1}(t)\cdot
\frac{|g^{(k)}(t)|}{1+|g(t)|^\alpha}\\
&=& \liminf_{t\to+\infty} \prod_{j=1}^{q-1}\log_j t\cdot t^k\cdot
\frac{f^{(k)}(e^t)\cdot e^{kt}+\sum_{j=1}^{k-1} c_j
f^{(j)}(e^t)\cdot e^{jt}}{1+|f(e^t)|^\alpha}\\
&=& \liminf_{x\to+\infty} \prod_{j=1}^{q-1}\log_{j+1} x\cdot (\log x)^k\cdot
\frac{f^{(k)}(x)\cdot x^k+\sum_{j=1}^{k-1} c_j
f^{(j)}(x)\cdot x^j}{1+|f(x)|^\alpha}\\
&\ge& \liminf_{x\to+\infty} \prod_{j=2}^{q}\log_j t\cdot \log x\cdot
\frac{f^{(k)}(x)\cdot x^k}{1+|f(x)|^\alpha}\\
&=& \liminf_{x\to+\infty} \frac{P_q(x)\cdot f^{(k)}(x)}{1+|f(x)|^\alpha},
\eeqaro
as desired. 

{\bf Case 2: } $L_{k-1}=\dots=L_{m+1}=0$, but $L_m\ne0$ for some
integer $m\ge0$, $m\le k-2$.

Then for $j=k-1,k-2,\dots,m+1$ and all large enough $x$ there is a
$\zeta_x\in[x,2x]$ such that
\beq\label{Case2a}
x\cdot |f^{(j)}(2x)| 
\le x\cdot |f^{(j)}(\zeta_x)|
=|f^{(j-1)}(2x)|-|f^{(j-1)}(x)| 
\le |f^{(j-1)}(x)|;
\eeq
here we have used that $|f^{(j-1)}|$ is decreasing (since $f^{(j-1)}$
is monotonic and $L_{j-1}=0$) and that $f^{(j-1)}(2x)$ and
$f^{(j-1)}(x)$ have the same sign. 

By induction we obtain for all $x$ large enough
\beq\label{Case2b}
x^{k-1}\cdot |f^{(k-1)}(2^{k-1-m}x)| \le \frac{1}{2^{(k-1-m)(k-2-m)/2}}\cdot x^m\cdot f^{(m)}(x)|.
\eeq

{\bf Case 2.1: } $L_m\ne \pm\infty$, i.e. $L_m\in\re$. 

Then for all $x$ large enough we have 
$$|f(x)|\ge \frac{x^m}{2m!}\cdot L_m,$$
hence 
\beq\label{Case2.1}
0\le \prod_{j=1}^{q} \log_j x\cdot \frac{x^m}{1+|f(x)|^\alpha} 
\le \prod_{j=1}^{q} \log_j x\cdot \frac{x^m}{1+\kl
  \frac{x^m}{2m!}\cdot L_m\kr^\alpha} 
\abb 0 \qquad (x\to\infty).
\eeq
From (\ref{MVT}) and (\ref{Case2b}) we conclude that for all $n$ large
enough
\beqaro
n^k\cdot |f^{(k)}(\zeta_n)|
&=& n^{k-1} |f^{(k-1)}(2n)-f^{(k-1)}(n)| \\
&\le & n^{k-1} |f^{(k-1)}(n)| \\
&=& 2^{(k-1-m)(k-1)}\cdot \kl\frac{n}{2^{k-1-m}}\kr^{k-1}
|f^{(k-1)}(n)| \\
&\le & 2^{(k-1-m)(k-1)}\cdot \kl\frac{n}{2^{k-1-m}}\kr^{k-1}
|f^{(m)}\kl \frac{n}{2^{k-1-m}}\kr|. 
\eeqaro
If we combine this estimate with (\ref{Case2.1}) and observe that
$f^{(m)}$ is bounded (since $L_m\in\re$), we obtain (with $C_m:=2^{(k-1-m)^2+k}$)
\beqaro
0 &\le & \prod_{j=1}^{q} \log_j \zeta_n\cdot \frac{\zeta_n^k \cdot  |f^{(k)}(\zeta_n)|}{1+|f(\zeta_n)|^\alpha} \\
&\le& \prod_{j=1}^{q} \log_j \zeta_n\cdot 2^k\cdot \frac{n^k \cdot  |f^{(k)}(\zeta_n)|}{1+|f(\zeta_n)|^\alpha} \\
&\le& C_m\cdot \prod_{j=1}^{q} \log_j \zeta_n\cdot \frac{n^m}{1+|f(\zeta_n)|^\alpha}\cdot  |f^{(m)}\kl\frac{n}{2^{k-1-m}}\kr|\\
&\le& C_m\cdot \prod_{j=1}^{q} \log_j \zeta_n\cdot\frac{\zeta_n^m}{1+|f(\zeta_n)|^\alpha} \cdot|f^{(m)}\kl\frac{n}{2^{k-1-m}}\kr|
\abb0 \quad (n\to\infty)
\eeqaro
for all $n$ large enough. This settles Case 2.1.

{\bf Case 2.2: } $L_m= \pm\infty$, w.l.o.g. $L_m=+\infty$. 

Then for all $x$ large enough we have $f^{(m)}(x)\ge
m!+1$, $f^{(m-1)}(x)\ge m!\cdot x+1,\dots, f'(x)\ge m\cdot x^{m-1}+1$ and finally
\beq\label{Case2.2a} 
f(x)\ge x^m,
\eeq
hence
$$\prod_{j=1}^{q} \log_j x\cdot \frac{x^m}{1+|f(x)|^\alpha}\abb0 \quad
(x\to\infty).$$
For $j=1,\dots,m$, by the Mean Value Theorem we find numbers $\zeta_x\in [x,2x]$
such that for all $x$ large enough
$$f^{(j-1)}(2x)=f^{(j-1)}(x)+x\cdot f^{(j)}(\zeta_x) \ge 0+ x\cdot f^{(j)}(x), $$
and by induction we conclude that
\beq\label{Case2.2b}
f(2^m x)\ge x^m\cdot f^{(m)}(x),
\eeq
provided that $x$ is large enough. On the other hand, $f^{(m+1)}$ is
positive and decreases to 0, so for a suitably chosen $x_0\ge0$ and all
$x\ge 2x_0$ we obtain
\beqaro
f^{(m)}(2^m x) 
&\le& f^{(m)}(x_0+2^m x) 
=f^{(m)}(x_0) +\int_{x_0}^{x_0+2^{m+1}\cdot\frac{x}{2}} f^{(m+1)}(t)\;dt\\
&\le& f^{(m)}(x_0) +2^{m+1}\cdot \int_{x_0}^{x_0+\frac{x}{2}}
f^{(m+1)}(t)\;dt\\
&=&2^{m+1}\cdot f^{(m)}\kl x_0+\frac{x}{2}\kr-(2^{m+1}-1)\cdot
f^{(m)}(x_0)\\
&\le& 2^{m+1}\cdot f^{(m)}(x)+0.
\eeqaro
Combining this with (\ref{Case2.2b}), we obtain for all $x$ large enough
$$2^{m+1}\cdot f(2^m x)\ge x^m \cdot f^{(m)}(2^m x),$$
hence (by replacing $2^m x$ with $x$)
\beq\label{Case2.2c} 
2^{m^2+m+1}\cdot f(x)\ge x^m \cdot f^{(m)}(x).
\eeq
If we combine this estimate with (\ref{MVT}), (\ref{Case2b}) and
(\ref{Case2.2a}), as in Case 2.1 we obtain
\beqaro
0&\le& 
P_q(\zeta_n)\cdot \frac{|f^{(k)}(\zeta_n)|}{1+|f(\zeta_n)|^\alpha}\\
&\le& C_m \cdot \prod_{j=1}^{q}\log_j \zeta_n \cdot\frac{n^m\cdot
  \l|f^{(m)}\kl\frac{n}{2^{k-1-m}}\kr\r|}{1+|f(\zeta_n)|^\alpha}\\
&\stackrel{(\ref{Case2.2c})}{\le}& C_m' \cdot \prod_{j=1}^{q}\log_j \zeta_n 
\cdot\frac{\l|f\kl\frac{n}{2^{k-1-m}}\kr\r|}{1+|f(\zeta_n)|^\alpha}\\
&\le& C_m' \cdot \prod_{j=1}^{q}\log_j \zeta_n \cdot |f(\zeta_n)|^{1-\alpha}\\
&\stackrel{(\ref{Case2.2a})}{\le}& C_m' \cdot \prod_{j=1}^{q}\log_j
\zeta_n \cdot \zeta_n^{m(1-\alpha)} 
\abb 0 \qquad (n\to\infty),
\eeqaro
where $C_m'$ is an appropriate constant. This settles this case as
well. 

{\bf Case 3:} $L_{k-1}=\dots=L_0=0$

In this case, (\ref{Case2b}) holds as well (with $m=1$), i.e.
$$|f'(x)|\ge x^{k-2}\cdot |f^{(k-1)}\kl 2^{k-2}x\kr|$$
for all $x$ large enough. Now we use 
$$|f^{(k)}(x)|\ge \frac{\varepsilon}{x^k \prod_{j=1}^{q} \log_j x}$$
(which is valid for all large enough $x$) and once more the Mean Value
Theorem to obtain for all large enough $x$ 
\beqaro
|f'(x)|&\ge& x^{k-2}\cdot|f^{(k-1)}(2^{k-2}x)-f^{(k-1)}(2^{k-1}x)|\\
&=& 2^{k-2}\cdot x^{k-1}\cdot |f^{(k)}(\zeta_x)| \qquad\qquad\qquad\qquad \mbox{(where
  $2^{k-2}x\le\zeta_x\le 2^{k-1} x$)}\\
&\ge& \frac{2^{k-2}\cdot x^{k-1}\cdot\varepsilon}{\zeta_x^k\cdot\prod_{j=1}^{q}\log_j \zeta_x}\\
&\ge& \frac{2^{k-2}\cdot  x^{k-1}\cdot\varepsilon}{(2^{k-1}x)^k\cdot\prod_{j=1}^{q}\log_j (2^{k-1} x)}\\
&\ge& c\cdot\frac{1}{x\cdot\prod_{j=1}^{q}\log_j x}
\eeqaro
with a suitable constant $c>0$, hence by integration
$$|f(x)|\ge c\cdot \log_{q+1} x+d\to\infty \qquad (x\to\infty),$$
(with some $d>0$), since $\frac{d}{dx}\log_{q+1} x=\frac{1}{x\prod_{j=1}^{q} \log_j
  x}$. This contradicts $L_0=0$, i.e. this case cannot
occur\footnote{In fact, Case 3 is the only part of the proof where it
  is crucial that in the assertion only the factors $\log_j x$ and not
  $(\log_j x)^\beta$ with $\beta>1$ occur. It would not work with
  $\beta>1$ since the improper integral
  $\int_{x_0}^{\infty}\frac{1}{x\log x\cdot\ldots\cdot \log_q x\cdot
    (\log_p x)^\beta}\;dx$ (with $x_0$ large enough) converges.}. 

This completes the proof of (\ref{MainAssertion1}). 

Now (\ref{MainAssertion2}) is an easy consequence from
(\ref{MainAssertion1}) and from Darboux' intermediate value theorem
for derivatives. Indeed, if there a $x_0$ such that $f^{(k)}(x)\ge 0$
for all $x\ge x_0$ or $f^{(k)}(x)\le 0$ for all $x\ge x_0$,
(\ref{MainAssertion2}) follows immediately from
(\ref{MainAssertion1}), applied to either $f$ or $-f$. Otherwise, by
Darboux's theorem there is a sequence $\gl x_n\gr_n$ tending to
$\infty$ such that $f^{(k)}(x_n)=0$ for all $n$, and
(\ref{MainAssertion2}) holds as well.  
\ebew

In view of Theorem \ref{mainresult} and the fact that the exponential
function grows larger than every polynomial, the following fact
certainly doesn't come as a big surprise: 

{\sl For every continuously differentiable function
  $g:[a,\infty)\abb\re$ we have 
\beq\label{Cor1}
\liminf_{x\to\infty} \frac{g'(x)}{e^{g(x)}}\le 0.
\eeq} 

Indeed, otherwise there would be an $\varepsilon>0$ and an $x_0\ge a$
  such that $g'(x)\ge \varepsilon \cdot e^{g(x)}$ for all $x\ge
  x_0$. In particular, $g'$ is positive on $[x_0,\infty)$, so $g$ is
  increasing there, hence $g'(x)\ge \varepsilon\cdot e^{g(x_0)}$ for
  all $x\ge x_0$, which implies $\lim_{x\to\infty} g(x)=\infty$. This
  enables us to conclude that
  $\frac{e^{g(x)}}{|g(x)|^2}\to\infty$ for
  $x\to\infty$. Combining this with the fact that
  $\liminf_{x\to\infty} \frac{g'(x)}{1+|g(x)|^2}\le 0$ by Theorem
  \ref{mainresult} gives the assertion.

However, it might be a bit surprising that this no longer holds if $g'$ is
replaced by higher derivatives of $g$, i.e. for $k\ge 2$ in general
the estimate $\liminf_{x\to\infty} \frac{g^{(k)}(x)}{e^{g(x)}}\le 0$
does not hold. This is demonstrated by the function $g(x):=-x^{k-3/2}$
which satisfies
$$\frac{g^{(k)}(x)}{e^{g(x)}} =
C\cdot \frac{x^{-3/2}}{\exp(-x^{k-3/2})}\abb\infty \qquad \mbox{ for } x\to\infty$$
with some $C>0$. 

On the other hand, for every $k$ times continuously differentiable function
  $g:[a,\infty)\abb\re$ ($k\ge 1$) we have
$$\liminf_{x\to\infty} \frac{g^{(k)}(x)}{1+e^{g(x)}}\le 0 
\qquad \mbox{ and } \qquad
\liminf_{x\to\infty} \frac{g^{(k)}(x)}{e^{|g(x)|}}\le 0.
$$
Both inequalities are proved by a similar reasoning as in the proof
  of (\ref{Cor1}), applying Theorem \ref{mainresult} with (for example)
  $\alpha=2$ and keeping in mind that $g^{(k)}(x)\ge \varepsilon$ for
  all $x\ge x_0$ would imply $g(x)\abb\infty$ for $x\to\infty$
  resp. that $x\mapsto \frac{e^{|g(x)|}}{1+|g(x)|^2}$ is bounded away
  from zero.  

\bibliographystyle{amsplain}

\begin{thebibliography}{99}
\bibitem{BarGrahlNevo} Bar, R.; Grahl, J.; Nevo, S.:  {\it
  Differential inequalities and quasi-normal families}, Anal. Math. Phys. {\bf 4} (2014), 63-71 

\bibitem{ChenNevoPang} Chen, Q.; Nevo, S.; Pang, X.-C.: {\it A general differential
  inequality of the $k$th derivative that leads to normality},
  Ann. Acad. Sci. Fenn. {\bf 38} (2013), 691-695

\bibitem{GrahlNevo-Spherical}
Grahl, J.; Nevo, S.: {\it Spherical derivatives and normal families},  
J. Anal. Math. {\bf 117} (2012), 119-128 

\bibitem{GN-Marty} Grahl, J.; Nevo, S.: {\it An extension of one direction
  in Marty's normality criterion}, Monatsh. Math. {\bf 174} (2014), 205-217 

\bibitem{GNP-NonExplicit} Grahl, J.; Nevo, S.; Pang, X.-C.: {\it A non-explicit
  counterexample to a problem of quasi-normality},
  J. Math. Anal. Appl. {\bf 406} (2013), 386-391 

\bibitem{LiuNevoPang} Liu, X.J., Nevo, S. and Pang, X.C.: {\it Differential
  inequalities, normality and quasi-normality}, Acta
  Math. Sin. (Engl. Ser.) {\bf 30} (2014), 277-282

\end{thebibliography}

\vspace{12pt}

\vspace{10pt}
\parbox{85mm}{\sl Jürgen Grahl\\
University of W\"urzburg \\
Department of Mathematics  \\     
W\"urzburg\\
Germany\\
e-mail: grahl@mathematik.uni-wuerzburg.de}
\hfill\parbox{75mm}{\sl Shahar Nevo \\
Bar-Ilan University\\
Department of Mathematics\\
Ramat-Gan 52900\\
Israel\\
e-mail: nevosh@math.biu.ac.il}

\end{document}